\documentclass{article}
\usepackage{amsmath,amsthm,amsfonts,amscd}

\numberwithin{equation}{section}

\def\cA{{\cal A}}
\def\cB{{\cal B}}
\def\cC{{\cal C}}
\def\cD{{\cal D}}

\def\cH{{\cal H}}
\def\cI{{\cal I}}

\def\cK{{\cal K}}
\def\cL{{\cal L}}
\def\cM{{\cal M}}

\def\cO{{\cal O}}

\def\cR{{\cal R}}
\def\cS{{\cal S}}

\def\cV{{\cal V}}

\def\bn{{\mathbb N}}
\def\br{{\mathbb R}}

\def\a{\alpha}
\def\b{\beta}
        
        \def\D{\Delta}
\def\eps{\varepsilon}

\def\k{\kappa}
\def\l{\lambda}       \def\L{\Lambda}
\def\m{\mu}

\def\x{\xi}
\def\p{\pi}

\def\s{\sigma}

\def\f{\varphi}
\def\F{\Phi}

\def\o{\omega}        \def\O{\Omega}

\def\gA{\mathfrak A}
\def\gB{\mathfrak B}
\def\gC{\mathfrak C}

\def\gS{\mathfrak S}
\def\imply{\Rightarrow}

\def\Ppo{{\cal P}_+^\uparrow}   
\def\sp{\mathrm{sp}}

\newtheorem{Thm}{Theorem}[section]
\newtheorem{Cor}[Thm]{Corollary}
\newtheorem{Prop}[Thm]{Proposition}
\newtheorem{Lemma}[Thm]{Lemma}
\newtheorem{Sublemma}[Thm]{Sublemma}
\theoremstyle{definition}

\theoremstyle{remark}
\newtheorem{rem}[Thm]{Remark} 

\renewcommand{\sectionmark}[1]{}

\begin{document}
\title{\huge Natural Energy Bounds\\
 in Quantum Thermodynamics}
\author{Daniele Guido$^{1)}$, Roberto Longo$^{2)}$}
\date{}
\markboth{}
{Energy bounds in quantum thermodynamics}
\maketitle
\bigskip\noindent
$^{1)}$ Dipartimento di Matematica, Universit\`a della Basilicata, I--85100
Potenza, Italy.\\
$^{2)}$ Dipartimento Matematica, Universit\`a di Roma ``Tor
Vergata'', I--00133 Roma, Italy.
\vskip 1cm

 \begin{abstract}
Given a stationary state for a noncommutative flow, we study a
boundedness condition, depending on a parameter $\b>0$,
which is weaker than the KMS equilibrium condition at inverse
temperature $\b$. This condition is equivalent
to a holomorphic property closely related to the one recently
considered by Ruelle and D'Antoni-Zsido and shared by a natural class of
non-equilibrium steady states.  Our holomorphic property is stronger
than the Ruelle's one and thus selects a restricted class of
non-equilibrium steady states.  We also introduce the complete
boundedness condition and show this notion to be equivalent to the
Pusz-Woronowicz complete passivity property, hence to the KMS
condition.

In Quantum Field Theory, the $\b$-boundedness condition can be
interpreted as the property that localized state vectors have energy
density levels increasing $\b$-subexponentially, a property which is
similar in the form and weaker in the spirit than the modular
compactness-nuclearity condition.  In particular, for a Poincar\'e
covariant net of C$^*$-algebras on the Minkowski spacetime, the
$\b$-boundedness property, $\b\geq 2\pi$, for the boosts
is shown to be equivalent to the Bisognano-Wichmann
property.  The Hawking temperature is thus minimal for a thermodynamical
system in the background of a Rindler black hole within the class of
$\b$-holomorphic states. More generally, concerning the Killing evolution
associated with a class of stationary
quantum black holes, we characterize KMS thermal equilibrium states
at Hawking temperature in terms of the boundedness property
and the existence of a translation symmetry on the horizon.

 \end{abstract}
 \vfill
\thanks{Work partially supported by MURST and GNAFA-INDAM}
 \newpage

 \setcounter{section}{-1}
\section{Introduction.}
\label{sec:intro}
In this paper we shall discuss a property for a state which is
invariant under a given one-parameter automorphism group of a
C$^*$-algebra. This property has two essentially
equivalent descriptions either as a boundedness condition or
as a holomorphic condition. The boundedness property has a natural
interpretation within Quantum Field Theory, being somehow similar to
the Haag-Swieca \cite{HS} compactness or Buchholz-Wichmann \cite{BW}
nuclearity-conditions, while the holomorphic condition
is a weakening of the KMS thermal equilibrium condition, related to
the conditions recently considered by Ruelle \cite{R} and D'Antoni-Zsido
\cite{DZ}, thus it naturally pertains to the context of Quantum Statistical
Mechanics.
It is then natural to discuss our property in a context where both the
two above subjects coexist: Black Hole Thermodynamics (see \cite{Wald2}).
Our main result is that for a translation invariant Quantum Field
Theory such property is equivalent to the KMS condition.

{\it Quantum field theory and the boundedness property.}
Let us consider a Quantum Field Theory on the Minkowski
spacetime and let $\gA(\cS)$ be the C$^*$-algebra
on the vacuum Hilbert space $\cH$ generated by
the observable localized in the region $\cS$. Clearly pure states on
the quasi-local observable algebra in the vacuum folium are given by
unit vectors of $\cH$, unique up to a phase, under the correspondence
\[
\xi\in\cH,\ ||\x||=1 \longrightarrow \o_\x
\]
where $\o_\x$ is the expectation functional $\o_\x(X)\equiv (X\x,\x).$
Denote by $\cL(\cS)$ the set of  vector states localized
in $\cS$, namely
\[
\cL(\cS)\equiv\{\xi\in\cH, ||\x||=1:
\o_\x |_{\gA(\cS')}=\o_0 |_{\gA(\cS')}\},
\]
with $\o_0=\o_\O$ the vacuum state.

It easy to see that $\xi\in\cL(\cS)$ if and only if there exists an
isometry $W\in\gA(\cS')'$ such that $\xi= W\O$; indeed $W$ is the
closure of the map $X'\O\to X'\x$, $X'\in\gA(\cS')$. Therefore, if
Haag duality holds for $\cS$, namely $\cA(\cS)=\gA(\cS')'$, where
$\cA(\cS)=\gA(\cS)''$ is the weak closure of $\gA(\cS)$, we have
\begin{equation}\label{isom}
\cL(\cS)=\{W\O: W\in\cA(\cS), W^*W =1\}\ .
\end{equation}
The $\b$-boundedness condition demands that all vectors in $\cL(\cS)$
have energy density levels increasing $\b$-subexponentially,
namely
\begin{equation}\label{bcond}
\int e^{-\lambda\b}\text{d}\mu_\x(\lambda) < +\infty,\ \forall \x\in
\cL(\cS),
\end{equation}
where $\mu_\x(B)=\int_B \text{d}(E(\lambda)\x,\x)$ is the spectral
measure associated by the Hamiltonian $H$ to $\xi$, $H=\int
\lambda\text{d}E(\lambda)$ being the spectral resolution of $H$. By the
spectral theorem, equation (\ref{bcond}) is equivalent to
$\cL(\cS)\subset D(e^{-\frac{\b}{2}H})$.

Since every operator $X\in\gA(\cS)$ with $||X||<1$ is a convex
combination of unitaries in $\gA(\cS)$ (see e.g \cite{Pe}), the
$\b$-boundedness condition is equivalent to $\gA(\cS)\O\subset
D(e^{-\frac{\b}{2}H})$ and turns out to be equivalent to
\begin{equation}
e^{-\b H}\gB_1\Omega\: \text{is a bounded set}
\label{mbound}
\end{equation}
for some, hence for all, weakly dense $^*$-subalgebra $\gB$ of
$\cA(\cS)$, where $\gB_1$ denotes the unit ball of $\gB$.

Equation (\ref{bcond}) already shows an interesting aspect of the
$\b$-boundedness property in Quantum Field Theory, namely that, as
mentioned, it can be formulated much in analogy with the
compactness-nuclearity \cite{HS,BW}.

Now, we do not have yet specified which region $\cS$ is supposed to
be.  It could be the entire Minkowski spacetime $M$, but an unbounded
region as a spacelike cone will already contain the relevant
information.  As $H$ commutes with all translations, property
(\ref{mbound}) for $\cS$ do implies the same property for the
translated regions $\{\cS+x,x\in\mathbb R^{d+1}\}$, hence for the
quasi-local algebra.  It is then almost immediate that the
$\b$-boundedness condition is equivalent to $H$ to be semibounded,
thus, by Poincar\'e covariance, to the positivity of the
energy-momentum.

Although instructive, the boundedness condition would certainly have a
limited interest if were confined to the above situation.  Such
condition acquires however a deeper role when the Hamiltonian is not
positive, as is the case in thermodynamical contexts.  To better
illustrate this point, we need to make a digression in
Statistical Mechanics and discuss a property equivalent to the
boundedness condition.

{\it $\b$-holomorphy and non-equilibrium states.}
As is well known, the thermal equilibrium states
in Quantum Statistical Mechanics at finite volume are the Gibbs states,
while at infinite volume
they are KMS \cite{HHW} states, the proper generalization of the
former, which are defined as follows.

Let $\gA$ be a C$^*$-algebra and $\alpha$ a one-parameter group of
automorphisms of $\gA$ and $\gB$ a dense $^*$-subalgebra of $\gA$.  A
state $\o$ of $\gA$, namely a normalized positive linear functional of
$\gA$, is a KMS state for $\a$ at inverse temperature $\b>0$ if, for
any $X,Y\in\gB$, the function
\begin{itemize}
\item[$(a)$] $F_{X,Y}(t)\equiv \o(\a_t(X)Y)$ extends to a function
in $A(S_\b)$,
\item[$(b)$] $F_{X,Y}(t+i\b)=\o(Y\a_t(X))$,
\end{itemize}
where $A(S_\b)$ is the algebra of functions analytic in the strip
$S_\b=\{0<\Im z <\b\}$, bounded and continuous on the closure $\bar
S_\b$.  The consideration of the subalgebra $\gB$ is indeed
unnecessary (but convenient for future reference) as properties $(a)$
and $(b)$ then hold for all $X,Y\in\gA$.

Let us now consider non-equilibrium statistical mechanics and in
particular the physical situation recently considered by Ruelle
\cite{R}.  There a quantum system $\Sigma$ is interacting with a set
of infinite reservoirs $R_{k}$ that are in equilibrium at different
temperatures $\b^{-1}_k$.  The system $\Sigma$ may be acted upon by a
force, that we assume to be time independent.  In this context, at
least if $\Sigma$ is finite, a natural class of stationary
non-equilibrium states occur, the non-equilibrium steady states.  If
we denote as above the observable C$^*$-algebra by $\gA$ and the time
evolution automorphism group by $\a$, a non-equilibrium steady state
$\o$ of $\gA$ satisfies property $(a)$ in the KMS condition, for all
$X,Y$ in a dense $^*$-subalgebra of $\gB$, with $\b=\text{min}\b_k$,
but not necessarily property $(b)$.

States with property $(a)$ have been independently also discussed by
D'Antoni and Zsido \cite{DZ}.  A typical example of a non-KMS state
satisfying property $(a)$ is provided by the tensor product of KMS
states at different temperatures; in this case the parameter $\b$ is
clearly the minimum of the inverse temperatures.  Further examples are
obtained by considering a KMS state with respect to a bounded
perturbation of the dynamics, showing a certain stability of the
holomorphic property.  In these examples the states satisfy property
$(a)$.

Because of property $(b)$, KMS states also satisfy the bound
\begin{itemize}
\item[$(c)$] $|F(t+i\b)|\leq C||X||\,||Y||$ for some constant $C>0$.
\end{itemize}
$X,Y\in\gB$, indeed one may take $C=1$.
We shall show that the bound $(c)$ automatically
occurs if property $(a)$ holds true for all the elements of a
C$^*$-algebra $\gA$.

We say that a $\a$-invariant state $\o$ is $\b$-holomorphic if
Properties $(a)$ and $(c)$ hold for all $X,Y$ in a dense
$^*$-subalgebra $\gB$ or, equivalently, if Property $(a)$ holds for
all $X,Y$ in the C$^*$-algebra $\gA$.

Since the bound $(c)$ does not necessarily hold for the
Non-Equilibrium Steady States considered by Ruelle \cite{R},
$\b$-holomorphic states form a subclass of the class of such
non-equilibrium states that, in a sense, are closer to equilibrium.

{\it KMS and complete $\b$-boundedness.} Because of the bound $(c)$,
the $\b$-holo\-mor\-phic property for a state $\o$ turns out to be
equivalent to the $\frac{\b}{2}$-boundedness property for the
Hamiltonian $H$ in the GNS representation.  The two properties are
thus two different aspects of the same notion.

Furthermore, if the constant $C$ is equal to $1$, we may use an
inequality of Pisier \cite{Pi1}, improved by Haagerup \cite{H}, to get
the inequality
\[
e^{\mp\b H}\leq 1+\D^{\pm 1}\ ,
\]
where $\D$ is the modular operator associated with the GNS vector $\O$
that, for simplicity, we are assuming to be separating for the weak
closure $\pi_\o(\gA)''$.  Thus the $\b$-holomorphic property with
$C=1$ entails that the Hamiltonian $H$ is dominated, in the above
sense, by the thermal equilibrium Hamiltonian $\log\D$.

A better understanding of the $\b$-holomorphic property is then
obtained by comparing it with the passivity condition of
Pusz-Woronowicz \cite{PW}, which is an expression of the second
principle of thermodynamics.  As is known, the passivity condition is
weaker than the KMS condition, while the complete passivity turns out
to be equivalent to the KMS property at some inverse temperature
(possibly $0$ or $+\infty$).  In analogy, we define the complete
$\b$-holomorphic property and show this property to be equivalent to
complete passivity, thus to the KMS condition.

We now return to Quantum Field Theory in a more general context.

{\it Black hole thermodynamics: minimality of the Hawking temperature.}
During the past thirty years, a theory of black hole thermodynamics
has been developed much in analogy with classical thermodynamics, see
\cite{Wald2}.  In this new context the thermodynamical functions acquire a
new meaning, for example the entropy is proportional to the area of
the black hole \cite{B}, yet a ``generalized second law of
thermodynamics'' holds.

As derived by Hawking \cite{H}, the temperature appearing in this formula
is a true physical temperature, in other words black holes do emit a thermal
radiation, provided quantum effects are taken into accounts.

This effect, or its closely related Unruh effect in the Rindler
spacetime \cite{U}, has been noticed \cite{Sew} to be essentially
equivalent to the Bisognano-Wichmann property in Quantum Field Theory.

Let now a Quantum Field Theory on the Minkowski spacetime $M$ be specified
by the algebras $\gA(\cO)$ of the observables localized in the regions
$\cO$. As is known the Rindler spacetime can be identified with a
wedge region $W$ of $M$, say $W=\{x\in {\mathbb R}^{d+1}: x_1>|x_0|\}$.
Thus $W$ represents the exterior of a Rindler black hole.
The pure Lorentz transformation in the $x_1$-direction on $M$ leave
$W$ globally invariant and thus give rise to a one-parameter
automorphism group $\a$ of the C$^*$-algebra $\gA(W)$.

The Bisognano-Wichmann theorem shows that, if $\gA$ is generated by a
Wightman field, then the restriction to $\gA(W)$ of the vacuum state
satisfies the KMS condition with respect to $\a$ at inverse
temperature $\b=2\pi$. One can then explain the Unruh effect on this
basis, see \cite{Sew}.

We shall show that the $\b$-boundedness property, $\b\geq 2\pi$, in the
above Quantum Field Theory context, actually implies the KMS property
at $\b=2\pi$, namely the state is at thermal equilibrium.

In particular, for a Poincar\'e covariant net of local observable
algebras on the Minkowski spacetime, we obtain a characterization of
the Bisognano-Wichmann property: there should exist a spacelike cone
$\cS$ contained in the wedge $W$ and a weakly dense $^*$-algebra $\gB$
of $\cA(\cS)=\gA(\cS)''$ such that
\begin{equation}\label{bR}
e^{-\pi K}\gB_1\Omega\quad \text{is a bounded set}
\end{equation}
of the underlying Hilbert space.  Here $\Omega$ is the vacuum vector
and $K$ is the generator of the boost unitary group corresponding to
$W$, namely the Killing Hamiltonian for the Rindler space $W$.  Note
that, if the Bisognano-Wichmann property and the split property
\cite{DL} hold, then also the modular compactness condition holds
true, namely $e^{-\lambda K}\gB_1\Omega$ is a compact set for all
$0<\lambda<\pi$ \cite{BDL1} ($e^{-2\pi K}$ is then the modular
operator associated with $(\cA(W),\Omega)$).  Thus the boundedness
condition (\ref{bR}) is very similar in the form to the modular
compactness-nuclearity condition \cite{BDL1,BDL2}.

Stating our result in the setting of Rindler spacetime, we obtain
the following: if a state $\o$ is $\b$-holomorphic, with respect to
the Killing evolution, with parameter $\b^{-1}$ less or equal to the
Hawking temperature, then $\o$ is indeed a KMS state at Hawking
temperature.

This minimality character of the Hawking temperature may have a
further physical interpretation, which is however limited by the fact
that only $\b$-holomorphic states appear in the context, see the
conclusion at the end of this paper.

We then extend our analysis to stationary stationary black holes
described by a globally hyperbolic spacetime with bifurcate Killing
horizon, making use of the net of observables localized on the
horizon, as in \cite{GLRV,L5}.  We show that, assuming the
$\b$-boundedness, the KMS condition is equivalent to the existence of
a translation symmetry on the horizon.

Our paper is organized in two sections.  In the first one we deal with
C$^*$-algebras and automorphisms or endomorphisms; we discuss there
the basic structure provided by the $\b$-holomorphic and
$\b$-boundedness condition, that we later apply in Quantum Field
Theory in Section 2.  In this second section we first discuss our
results in the Minkowski, or Rindler, spacetime; then we study the
corresponding structure for one-dimensional nets and we then apply
these results to the case of a Quantum Field Theory on a spacetime
with bifurcate Killing horizon.

\section{Holomorphic states and the boundedness property}
In this section we discuss the basic structure of the holomorphic and
of the boundedness properties.  In the next section we shall apply our
results to the Quantum Field Theory context.
\subsection{General properties}
Let $\cH$ be a Hilbert space, $\O\in\cH$ a vector and $K$ a
selfadjoint operator on $\cH$.

We shall say that a linear subspace $\gS$ of $B(\cH)$ is {\it
$\b$-bounded} with respect to $K$ and $\O$ if $\gS\O\subset D(e^{-\b
K})$, the domain of $e^{-\b K}$, and the linear map
$$\Phi_{\b}: X\in\gS\to e^{-\b K}X\O\in\cH$$
is bounded.

Note that if $\gS$ is $\b$-bounded  then it is $\b'$-bounded for
all $0<\b'\leq\b$, indeed if $E$ is the spectral projection of $K$
relative to the interval $[0,+\infty)$, then
\begin{multline}\label{ub}
||e^{-\b' K}X\O||^2
=||Ee^{-\b' K}X\O||^2+||(1-E)e^{-\b' K}X\O||^2\\
\leq ||X\O||^2 + ||(1-E)e^{-\b K}X\O||^2
\leq ||X||^2 + ||e^{-\b K}X\O||^2
\leq 1+\|\Phi_{\b}\|^2,
\end{multline}
for all $X\in\gS_1$, where $\gS_1$ denotes the unit ball of $\gS$.

\begin{Lemma}\label{closable} Let $\gS$ be a $\b$-bounded linear subspace of
$B(\cH)$ for some $\b>0$.  Then
the restriction of $\Phi=\Phi_{\b}$ to $\gS_1$ is continuous with the
weak-operator topology on $\gS$ and the weak topology on $\cH$.

If $\gS$ is norm closed and $\gS\O$ is contained in $D(e^{-\b K})$,
then $\gS$ is automatically $\b$-bounded.
\end{Lemma}
\begin{proof} Assume then that $\Phi$ is bounded.  To check the
continuity, let $X_i\in\gS_1$ be a net weakly convergent to $0$.  Then
for all $\xi\in D(e^{-\b K})$ we have
\[
\lim_i (\xi,\Phi(X_i))=\lim_i (\xi,e^{-\b K}X_i\O)
=\lim_i (e^{-\b K}\xi,X_i\O)=0\ .
\]
as $D(e^{-\b K})$ is dense in $\cH$ and $\{\Phi(X_i)\}_i$ is bounded
by assumption, it follows that $\Phi(X_i)\to 0$ weakly.

It remains to show that $\Phi$ is bounded if $\gS$ is norm closed.
We shall show that $\Phi$ is closable, thus $\Phi$ will
be bounded by the closed graph theorem.

Let $\{X_n\}$ be a sequence in $\gS$ and $\eta\in\cH$ be such that
$X_n\to 0$ and $\Phi(X_n)\to \eta$ in norm. As $\{X_n\}$ is a bounded
set, by the just proved weak continuity of $\Phi$ we have that
$\Phi(X_n)\to 0$ weakly, thus $\eta = 0$.
\end{proof}

\begin{Lemma}\label{1.1} Let $\cM$ be a von Neumann algebra on the
Hilbert space $\cH$, $\O$ a vector, $K$ a selfadjoint operator of
$\cH$ and $\b>0$.

The following are equivalent:
\begin{itemize}
\item[$(i)$] There exists a weakly dense $^*$-subalgebra $\gB$ of
$\cM$  which is $\b$-bounded (w.r.t. $K$ and $\O$);
\item[$(ii)$] $\cM$ is $\b$-bounded;
\item[$(iii)$] $\gA\O\subset D(e^{-\b K})$, where $\gA$ is some weakly
dense C$^*$-subalgebra of $\cM$.
\end{itemize}
In this case $\F_{\b}|_{\gB}$ and $\F_{\b}|_{\cM}$ have the same norm.
\end{Lemma}
\begin{proof} Assuming in $(i)$ that $e^{-\b K}\gB_1\O$ is contained in
the ball of radius $C>0$, we shall show that the same is true for
$e^{-\b K}\cM_1\O$, i.e. $(ii)$ holds.
Let $X\in\gB$. By Kaplanski density theorem \cite{Takesaki}, there
exists a net  of operators $X_i\in\gB_1$ strongly convergent to $X$.
Since $\|e^{-\b K}X_i\O\|\leq C$, we
may assume, possibly restricting to a subnet, that $e^{-\b K}X_i\O$ weakly
converges to $\eta\in\cH$, $\|\eta\|\leq C$. Now take
$\xi\in\cD(e^{-\b K})$. We have
$$
(\xi,\eta)=\lim_i(\xi,e^{-\b K}X_i\O)=
\lim_i(e^{-\b K}\xi,X_i\O)=(e^{-\b K}\xi,X\O)\ .
$$
Since $e^{-\b K}$ is self-adjoint, this means that $X\O\in\cD(e^{-\b K})$
and $e^{-\b K}X\O=\eta$, i.e. $\|e^{-\b K}\cM_1\O\|\leq C$.

Now $(ii) \Rightarrow (iii)$ and $(iii) \Rightarrow (i)$ follows by
Lemma \ref{closable} as $\gA$ is norm closed.
\end{proof}

Let now $\gA$ be a unital C$^*$-algebra, $\a$ a one-parameter
automorphism group of $\gA$ and $\o$ a $\a$-invariant state of $\gA$.
We shall always assume that the maps $t\in \mathbb R\to\o(\a_t(X)Y)$
are continuous for all $X,Y\in\gA$.
Denote by $(\cH,\pi,\O)$ the GNS triple associated with $\o$ and by
$U$ the one-parameter unitary group on $\cH$ implementing $\pi\cdot\a$:
\[
U(t)\pi(X)\O=\pi(\a_t(X))\O\ ,\quad X\in\gA\ .
\]
>From our continuity assumption, it follows at once that $U$ is
strongly continuous.
Denote by $K$ the infinitesimal generator of $U$.
\begin{Prop}\label{pure} If $\o$ is a pure state, $\p(\gA)$ is
$\b$-bounded with respect to $K$ and $\O$ if and only if
the spectrum of $K$ is bounded below.
\end{Prop}
\begin{proof} If $\o$ is pure, then $\cM\equiv\pi(\gA)''=B(\cH)$,
hence, if the boundedness condition holds,
$e^{-\b K}B(\cH)_1\O$ is bounded. As $B(\cH)_1\O=\cH_1$, it follows
that $e^{-\b K}$ is bounded, thus $K$ is semibounded. The converse is
obvious.
\end{proof}

Let $\gB$ a $^*$-subalgebra of $\gA$.
The state $\o$ of $\gA$ is {\it $\b$-holomorphic} on $\gB$ if $\o$ is
$\a$-invariant and for every $X,Y\in\gB$ the function
$F_{X,Y}(t)=\o(\a_t(X)Y)$ is the boundary value of a function in
holomorphic in the strip $S_\b=\{z:0<\Im z <\b\}$, continuous in $\bar S_\b$.

Denote by $A(S_\b)$ the algebra of functions holomorphic in
$S_\b$, bounded and continuous in its closure $\bar S_\b$.
\begin{Prop}\label{ab}
Let $\gA$, $\a$, $\o$, $K$, $\O$ as before,
$\gB$ a $^{*}$-subalgebra of $\gA$.  Then $\o$ is
$\b$-holomorphic on $\gB$ if and only if $\p(\gB)\O\subset
D(e^{-\frac{\b}{2}K})$.

In this case $F_{X,Y}$ extends to a function in $A(S_\b)$, for all
$X,Y\in\gB$.
\end{Prop}
\begin{proof} Immediate by Lemma \ref{aLemma}.
\end{proof}
\begin{Thm}\label{equiv} Let $\gA$, $\a$, $\o$, $K$, $\O$ as before,
$\gB$ a $^{*}$-subalgebra of $\gA$.  Then $\p(\gB)$ is
$\frac{\b}{2}$-bounded w.r.t. $K$ and $\O$ if and only if $\o$ is
$\b$-holomorphic on $\gB$ and
\begin{equation}\label{bddholo}
|F_{X,Y}(t+i\b)|\leq C \|X\|\ \|Y\|,
\end{equation}
for some constant $C>0$.

If moreover $\gB$ is norm closed, then $\pi(\gB)$ is
$\b/2$-bounded w.r.t. $K$ and $\O$ if and only if $\o$ is
$\b$-holomorphic for $\gB$.
\end{Thm}

\begin{proof} Clearly $F_{Y^*,X}(t)=(e^{itK}\p(X)\O,\p(Y)\O)$ for all
$X,Y\in\gA$.  If $\o$ is $\b$-holomorphic on $\gB$ and $X,Y\in\gB$, then
$F_{Y^*,X}$ is boundary value of a function in $A(S_{\b})$,
thus by Lemma \ref{aLemma} $\pi(\gB)\O\subset D(e^{-\frac{\b}{2}K})$.
If moreover  the bound  (\ref{bddholo}) holds, then
\begin{equation*}
\|\F_{\b/2}|_{\gB}\|^{2}
=\sup_{X,Y\in\gB_{1}}|(e^{-\frac{\b}{2}K}\p(X)\O,e^{-\frac{\b}{2}K}\p(Y)\O)|=
\sup_{X,Y\in\gB_{1}}|F_{Y^*,X}(i\b)|\leq C\ ,
\end{equation*}
thus $\gB$ is $\b/2$-bounded.

Conversely, if $\gB$ is $\b/2$-bounded,
then $\pi(\gB)\O\subset D(e^{-\frac{\b}{2}K})$ and the same
computation done above yields, by Lemma \ref{aLemma},
that $F_{X,Y}\in A(S_\b)$ for all $X,Y\in\gB$
and
$$
\sup_{X,Y\in\gB_{1}}|F_{X,Y}(t+i\b)|
=\sup_{X,Y\in\gB_{1}}|F_{X,Y}(i\b)|=\|\F_{\b/2}|_{\gB}\|^{2}\ .
$$
so the the bound (\ref{bddholo}) holds with $C=\| \F_{\b/2}|_{\gB}\|^2$.

Since by Prop.  \ref{ab} $\o$ is $\b$-holomorphic on $\gB$ iff
$\pi(\gB)\O\subset D(e^{-\frac{\b}{2}K})$, the
rest follows by Lemma \ref{1.1}.
\end{proof}
Of course KMS states at inverse temperature $\b'>0$ are
$\b$-holomorphic for all $0<\b\leq\b'$ on all $\gA$ and satisfy the
bound (\ref{bddholo}).
\begin{Cor} With the above notations, let $\o$ be $\b$-holomorphic
on $\gB$. Then $\o$ is $\b$-holomorphic
on the closure $\bar\gB$ iff the bound $(\ref{bddholo})$ holds.
\end{Cor}
\begin{proof} Immediate by  Lemma \ref{1.1} and Theorem \ref{equiv}.
\end{proof}
\subsection{Complete $\beta$-holomorphy and KMS condition}
 We begin to recall the following inequality (\ref{PH}) which is due
 to Pisier, with the improved constant due to Haagerup, see
 \cite{Pi}.

 \begin{Thm}\label{Pisier}{\rm \cite{Pi1,Haa}.} Let $\Phi$
 be a bounded linear map from a C$^*$-algebra $\gA$ to a Hilbert space
 $\cH$.  Then there exist two states $\f$ and $\psi$ on $\gA$ such
 that
\begin{equation}\label{PH}
	\|\Phi(X)\|^2\leq
\|\Phi\|^2(\f(X^* X)+\psi(XX^*)),\quad X\in\gA\ .
\end{equation}
In the special case where $\gA$ is unital
and $||\Phi||=||\Phi(1)||$, one may take
$\f=\psi=||\Phi||^{-1}(\Phi(\cdot),\Phi(1))$ in eq. (\ref{PH}).
\end{Thm}
\noindent The special case in last part of the statement is obtained
during Haagerup's proof of the inequality (\ref{PH}); such proof is
not difficult and can be found in \cite{Pi}, Thm 7.3.
\begin{Cor}\label{ineq} Let $\cM$ be a von Neumann algebra on a
Hilbert space $\cH$, $\O\in\cH$ a cyclic unit vector
for $\cM$ and $U(t)=e^{itK}$ a $\O$-fixing one-parameter unitary group
on $\cH$ implementing automorphisms of $\cM$.

If $||e^{-\b K}\cM_1\O||\leq 1$, then
\begin{equation}\label{leq}
e^{-2\b K}\leq 1+\D E\ ,
\end{equation}
where $E\in\cM$ is the projection onto $\cH_0\equiv\overline{\cM'\O}$ and
$\D$ is the modular operator on $\cH_0$ associated with $(E\cM E,\O)$.
\end{Cor}
\begin{proof}
As the map $\Phi_\b:X\in\cM\to e^{-\b K}X\O\in\cH$ satisfies
$||\Phi_\b||=||\Phi_\b(1)||=1$, the inequality (\ref{PH}) holds with
$\f=\psi=\o$ (see \ref{Pisier}),
where $$\o(X)=(\Phi_\b(X),\Phi_\b(1))=(X\O,\O),$$ namely
\begin{equation}\label{dom}
||e^{-\b K}X\O||^2\leq ||X\O||^2 + ||X^*\O||^2\ , X\in\cM\ .
\end{equation}
Assuming first that $\O$ is also separating, i.e.  $E=1$ ,if $X\in\cM$
and $X\O\in D(\D)$ we then have
\[
(e^{-\b K}X\O,e^{-\b K}X\O)\leq (X\O,X\O) + (\D X\O,X\O)=((1+\D)X\O,X\O)\ .
\]
Since $U$ implements automorphims of $\cM$ and $U(t)\O=\O$,
by the modular theory $U(t)$ and $\D^{is}$ commute, thus there exists
a strongly dense subalgebra $\gB$ of $\cM$ such that $\gB\O$ is a
core for every continuous function of $K$ or of $\log\D$.
Taking $X\in\gB$, the above equation gives
\[
(e^{-2\b K}X\O,X\O)\leq ((1+\D)X\O,X\O)\ ,
\]
thus $e^{-2\b K}\leq 1+\D$ as $\gB\O$ is a core for both
$e^{-2\b K}$ and $1+\D$.

In general, if $E\neq 1$, we may consider the the reduced von Neumann
algebra $E\cM E$ on $\cH_0$.  Since $\O$ is separating for $E\cM E$
and $K$ commutes with $E$, the above shows that
\begin{equation}\label{fi}
e^{-2\b K}E\leq E+\D E \ .
\end{equation}
Thus, by the eq. (\ref{dom}), we have for all $X\in\cM$
\begin{equation*}
||e^{-\b K}(1-E)X\O||^2\leq ||(1-E)X\O||^2 + ||X^*(1-E)\O||^2
= ||(1-E)X\O||^2\ ,
\end{equation*}
that entails
\begin{equation}\label{si}
e^{-2\b K}(1-E)\leq 1-E \ .
\end{equation}
Combining the inequalities (\ref{fi}) and (\ref{si}) we get the
desired inequality (\ref{leq}).
\end{proof}

\begin{rem}
In general, the tensor product of bounded maps from C$^{*}$-algebras
to Hilbert spaces is not bounded.  Indeed, if $\cM$ is a von Neumann
algebra with a cyclic and separating vector $\O$, the map $\Phi_{\a}:
X\in\cM\mapsto \D^{\a}X\O$ is bounded for any $\a\in[0,1/2]$, but
$\Phi_{\a}\otimes\Phi_{\b}$ is not necessarily bounded, if $\a\ne\b$.
	
As an example, let $\cM\simeq B(\cH)$ be a type I$_{\infty}$ factor
acting by right multiplication on the Hilbert space $L^{2}(\cH)$ of
the Hilbert-Schmidt operators affiliated on $\cH$, $h$ a positive
Hilbert-Schmidt operator of norm 1 in $L^{2}(\cH)$.  Then $\D
X=h^{2}Xh^{-2}$, $X\in L^{2}(M)$, and
$\Phi_{\a}(X)=h^{2\a}Xh^{1-2\a}$, $X\in \cM$, $X\in L^{2}(\cH)$.  As a
consequence, if $h^{2\a}$ and $h^{1-2\a}$ are both are Hilbert-Schmidt,
the norm of the map $\Phi_{\a}$ is bounded by
$\|h^{2\a}\|_{2}\cdot\|h^{1-2\a}\|_{2}$.  Therefore,
$||\Phi_{\a}\otimes\Phi_{\b}||\leq\|h^{2\a}\|_{2} \cdot
\|h^{1-2\a}\|_{2} \cdot \|h^{2\b}\|_{2} \cdot \|h^{1-2\b}\|_{2}$.  For
example, if the eigenvalue sequence of $h$ is $\{2^{-n}\}_{n\in\bn}$,
$\Phi_{\a}\otimes\Phi_{\b}$ is bounded for every $\a,\b\in(0,1/2)$.
	
On the other hand, denoting by $F$ the unitary operator in
$\cM\otimes \cM$ that, under the isomorphism with
$B(\cH\otimes\cH)$, is defined by $F\xi\otimes \eta=\eta\otimes \xi$,
one has $F(X\otimes Y)=(Y\otimes X)F$, hence
$$
\Phi_{\a}\otimes\Phi_{\b}(F) =h^{2\a}\otimes h^{2\b}Fh^{1-2\a}\otimes
h^{1-2\b} =h^{1-\eps}\otimes h^{1+\eps}F
$$
where $\eps=2\a-2\b$.  Thus, if $h^{1-\eps}$ or $h^{1+\eps}$ are not
Hilbert-Schmidt, $\Phi_{\a}\otimes\Phi_{\b}$ is unbounded.  For
example, if the eigenvalue sequence of $h$ is $\{\frac{1}{\sqrt n\log
n}\}_{n\in\bn}$, then $\Phi_{\a}\otimes\Phi_{\b}$ is bounded if and only if
$\a=\b$.
\end{rem}

Note that, if in Corollary \ref{ineq} the vector $\O$ is separating
for $\cM$, then by modular theory $J\D J =\D^{-1}$ and $J K J = -K$,
where $J$ is the modular conjugation of $(\cM,\O)$, thus the
inequality $e^{2\b K}\leq 1 +\D^{-1}$ holds too, therefore
\[
e^{\mp 2\b K}\leq 1 + \D^{\pm 1} \ ,
\]
showing that the $\b$-boundedness condition with constant $C=1$ sets a
a bound on the Hamiltonian $K$ by the equilibrium Hamiltonian $\log\D$.

It is then natural to look for further conditions that entail $K$
to be proportional to $\log\D$.

Let $\gA$ be a C$^*$-algebra acting on a Hilbert space $\cH$ with a
cyclic vector $\O\in\cH$ and $U(t)=e^{itK}$ a $\O$-fixing
one-parameter unitary group on $\cH$ implementing automorphisms of
$\gA$.  We shall say that $\gA$ is \emph{completely $\b$-bounded} with
respect to $K$ and $\O$ if
$||\Phi_\b\otimes\Phi_\b\otimes\cdots\otimes\Phi_\b||\leq 1$, for all
finitely many tensor products of $\Phi_\b$ with itself, with
$\Phi_\b:X\in\gA\to e^{-\b K}X\O\in\cH$.

Here we consider the spatial tensor product norm on
$\gA\otimes\gA\otimes\cdots\otimes\gA$ and, of course, the Hilbert
space tensor norm on $\cH\otimes\cH\otimes\cdots\otimes\cH$.

Note that, by Lemma \ref{1.1}, the complete $\b$-boundedness
condition can be equivalently formulated in terms of the von Neumann
algebra $\cM=\gA''$.

\begin{Cor}\label{T}
Let $\cM$ be a von Neumannalgebra acting cyclic and separating vector
$\O\in\cH$ and $U(t)=e^{itK}$ a $\O$-fixing one-parameter unitary
group implementing automorphisms of $\gA$.  If $\cM$ is completely
$\b$-bounded with respect to $K$ and $\O$, then
\[
2\b K= -T\log\D\ ,
\]
where $T$ is a positive linear operator, $0\leq T\leq 1$, commuting
with $\D$, $K$ and the modular conjugation $J$ of $(\cM,\O)$.
\end{Cor}
\begin{proof}
By Corollary \ref{ineq} applied to the $n$-fold tensor product,
we have
$e^{-2\b K}\otimes\cdots\otimes e^{-2\b K}\leq
1+\D\otimes\cdots\otimes\D$. By the modular theory $K$ and $\D$ commute,
thus a simple application of the Gelfand-Naimark theorem implies
$e^{-2n\b K}\leq 1+\D^n$
for all $n$, thus $e^{-2\b K}\leq \sqrt[n]{1+\D^n}$.
Taking the limit as
$n\to\infty$ we obtain
\[\label{leq2}
e^{-2\b K}\leq \text{max}(1,\D).
\]
As $JKJ=-K$ and $J\log\D J= -\log\D$, the above inequality also
entails that
\[
e^{2\b K}\leq \text{max}(1,\D^{-1}).
\]
Taking logarithms, as $K$ and $\D$
commute, these inequalities
imply respectively $2\b K E_+ \leq -\log\D  E_+$ and
$2\b K E_- \geq -\log\D  E_-$, where $E_{+/-}$ are the spectral
projections of $\log\D$ corresponding to the positive/negative
half-line.

The projection $E_0$ onto the kernel of $\log\D$ clearly commutes with
$\D$ and $K$, thus denoting by $T$ the closure of
$-2\b K(\log\D)^{-1}(1-E_0)$, $T$ is a bounded positive linear
contraction, $0\leq T\leq 1$, commuting both with $\D$ and $K$.
As $JKJ=-K$ and $J\log\D J= -\log\D$ and $JE_0 J=E_0$, we also have
$JTJ=T$.
\end{proof}
We now recall a weak form of the characterization of the KMS property in
terms of the Roepstoff-Araki-Sewell auto-correlation lower bound and
the Pusz-Woronowicz passivity condition, see \cite{BR}.
\begin{Thm}\label{PW}{\rm{\cite{AS,Ro,PW}}}.
Let $\cM$ be a von Neumann algebra and $\O$ a cyclic and separating
vector and $U(t)=e^{itK}$ a $\O$-fixing one-parameter unitary group
implementing  automorphisms of $\cM$. Consider following properties
\begin{itemize}
\item[$(i)$] $(K X\O,X\O)\geq 0$ for all
$X\in \cM_{sa}$ in the domain of the derivation $[\cdot,K]$, where
$\cM_{sa}$ denotes the selfadjoint real subspace of $\cM$.
\item[$(ii)$] $K = -\lambda\log\D$ for some $\lambda\in [0,+\infty)$.
\end{itemize}
Then $(ii)\Rightarrow (i)$. If the group of automorphisms of $\cM$
preserving $(\cdot \O,\O)$ is ergodic then $(i)\Rightarrow (ii)$.
\end{Thm}
Note that case $(i)$ in the above theorem cannot hold with $K$ a
non-trivial positive operator (ground state) because then $K$ is
affiliated to $\cM$ by a theorem of Borchers \cite{Borc2} (see also
\cite{L}), thus $K=0$ because $\O$ is separating.

We shall need to test the above inequality in $(i)$ for more vectors.
\begin{Lemma}
Let $\cM$ be a von Neumann algebra and $\O$ a cyclic and separating
vector. Then
\begin{equation}\label{pass}
-(\log \D \x,\x)\geq 0
\end{equation}
for all vectors $\x\in \cK\cap D(\log\D)$, where $\cK$ is the real
Hilbert subspace given by
$\cK\equiv\overline{\cM_{sa}\O}$.
\end{Lemma}
\begin{proof}
Let $S=J\D^{\frac{1}{2}}$ be the Tomita operator and
$E_n$ be the spectral projection of $\log\D$ corresponding to the interval
$(-n,n)$. Then $E_n$ commutes with the (real, unbounded)  projection
$P=\frac{1}{2}(1+S)$ onto $\cK$ (see Lemma \ref{TT}) because it commutes with
$S=J\D^{\frac{1}{2}}$: indeed $E_n$
commutes both with $\D$ and with $J$ (being a real even function
of $\log\D$).

Denoting by $\cM(-n,n)$ the space of elements of $\cM$
whose spectrum under the modular group $\sigma_t =$ Ad$\D^{it}$
lies in $(-n,n)$ and by $\cM_{sa}(-n,n)$ the real subspace of
selfadjoint elements of $\cM(-n,n)$, we have
\begin{equation}\label{P}
\cM_{sa}(-n,n)\O= \frac{1}{2}(1+S)\cM(-n,n)\O=P\cM(-n,n)\O\ .
\end{equation}
Let $\cH_S$ denote the Hilbert space $D(S)$ equipped with the $S$-graph
scalar product
\[
(\x,\eta)_S\equiv (\x,\eta) + (S\x,S\eta) =
(\x,\eta) + (\D^{\frac{1}{2}}\x,\D^{\frac{1}{2}}\eta),\ \x,\eta\in
D(S),
\]
and notice that any $\gS$  subset of $E_n\cH$ is contained in $\cH_S$
and its closure  $\overline{\gS}$ in $\cH$ coincides with its closure
in $\cH_S$ (because the restriction of $\D^{\frac{1}{2}}$ to
$E_n\cH$  is bounded).

Notice also that, as a linear operator of $\cH_S$, $P$ is bounded,
indeed $P$ is the (real) orthogonal projection of $\cH_S$ onto $\cK$.

Therefore, on the Hilbert space $\cH_S$, we have by eq. (\ref{P})
that
\[
\overline{\cM_{sa}(-n,n)\O}=P\overline{\cM(-n,n)\O}
=PE_n\cH_S=E_nP\cH_S= E_n\cK\ .
\]
As the the inequality (\ref{pass}) holds for all $\x\in
\cM_{sa}(-n,n)\O$, it then holds for all $\x\in E_n\cK$.

Given $\x\in\cK\cap D(\log\D)$ the sequence of vectors
$\x_n\equiv E_n\x$ belongs to $\overline{\cM_{sa}(-n,n)\O}=E_n\cK$ and
\begin{equation*}
||\x_n-\x||\to 0, \quad
||\log\D\x_n-\log\D\x||\to 0,
\end{equation*}
therefore
\[
-(\log \D \x,\x)= -\lim_n(\log \D \x_n,\x_n)\geq 0\ .
\]
\end{proof}

\begin{Lemma}\label{TT} Let $\cM$ be a von Neumann algebra with a cyclic
separating vector $\O$ and let $S$ be the associated
Tomita's operator, i.e. the closure of $X\O\to X^*\O$, $X\in\cM$.
Then $\overline{\cM_{sa}\O}=\{\x\in D(S): S\x=\x\}$.
\end{Lemma}

\begin{proof} As is well known, $S^*=F$, where $F$ is the
Tomita's operator associated with $\cM'$ and $\O$. With
$\cK\equiv\overline{\cM_{sa}\O}$, let $\tilde S$ the operator given by
$\tilde S:\x+ i\eta\to\x- i\eta$, $\x,\eta\in\overline{\cM_{sa}\O}$,
and define analogously $\tilde F$ with respect to $\overline{{\cM'}_{sa}\O}$.
Then $\tilde S\supset S$ and $\tilde F\supset F$ and
$\tilde S^*\supset \tilde F$, thus $\tilde S =S$ and $\tilde F =F$.
\end{proof}

For completeness we give a generalization of the inequality
(\ref{pass}) that, at the same time, gives a direct proof of it.
\begin{Prop}
Let $\cH$ be a complex Hilbert space and $\cK$ a standard real Hilbert
subspace, namely $\cK\cap i\cK=\{ 0 \}$, $\overline{\cK+i\cK}=\cH$. Then
\begin{equation}\label{pass2}
-(\log \D \x,\x)\geq 0
\end{equation}
for all vectors $\x\in \cK\cap D(\log\D)$, where $\D$ is the modular
operator on $\cH$ associated with $\cK$.
\end{Prop}
\begin{proof} Since the kernel of $\log\D$ is invariant under $J$ and
$\D$, one can decompose $\cK$ as a direct sum of two components, one
corresponding to the kernel of $\log \D$, and one to its orthogonal
complement, thus the inequality (\ref{pass2}) can be proved for each
component separately.  Since the inequality is obviously satisfied on
the kernel of $\log\D$, we may just suppose the kernel of $\log\D$ to
be trivial.

Now we give an explicit description of the vectors of $\cK$ (cf.
\cite{FG2}) which allows an immediate verification of the inequality
\ref{pass2}.

Let us chse choose a selfadjoint antiunitary $C$ commuting with $J$
and $\D$, and set $U=JC$, so that $U\log\D U=-\log\D$. Then denote
with $\cL$ the real vector space of $C$-invariant vectors in the
spectral subspace $\{\log\D>0\}$ and by $\psi^{\pm}$ the maps $\psi^{+}
:y\in \cL\mapsto U\cos\Theta y+\sin\Theta y$, $\psi^{-}:y\in\cL\mapsto
iU\cos\Theta y-i\sin\Theta y$, where the operator $\Theta$ is defined by
$|\log\D|=-2\log\tan\Theta/2$, $\s(\Theta)\subseteq[0,\pi/2]$.

Since $U$ maps the spectral space $\{\log\D>0\}$ onto the spectral
space $\{\log\D<0\}$, both $\psi^{+}$ and $\psi^{-}$ are isometries,
and a simple calculation shows that their ranges are real-orthogonal.
Moreover, decomposing $\cH$ as $\{\log\D<0\}\oplus\{\log\D>0\}$, one
can show that any solution of the equation $Sx=x$ can be written as a
sum $\psi^{+}(y)+\psi^{-}(z)$, namely the map $\psi^{-}+\psi^{+} :
\cL\oplus_{\br} \cL\to\cK$ is an isometric isomorphism of real Hilbert
spaces.

Moreover, for any $y,z\in \cL$, $(\psi^{\mp}(y),\log\D\psi^{\pm}(z))$ is
purely imaginary, therefore
\begin{multline}
((\psi^{-}(y)+\psi^{+}(z)),\log\D(\psi^{-}(y)+\psi^{+}(z)))\\
=(\psi^{-}(y),\log\D\psi^{-}(y))+(\psi^{+}(z),\log\D\psi^{+}(z)),
\end{multline}
namely the inequality should be checked on $\psi^{+}(L)$ and
$\psi^{-}(L)$ separately.  Finally, $$(\psi^{-}(y),\log\D\psi^{-}(y))
= -(y,\cos\Theta\log\D y)\leq0,$$ since $\cos\Theta$ and $\log\D$ are
commuting positive operators on $\cL$, and the same holds on the range
of $\psi^{+}$.
\end{proof}

\begin{Thm}\label{cb2}
Let $\cM$ be a von Neumann algebra with a cyclic and separating vector
$\O$ and $U(t)=e^{itK}$ a $\O$-fixing one-parameter unitary group
implementing automorphisms of $\cM$.  If $\cM$ is completely
$\b$-bounded with respect to $K$ and $\O$, then either $K=0$ or {\rm
Ad}$U(t)$ satisfies the KMS condition at some inverse temperature
$\b_0\geq\b$; indeed $\b_0$ is the greatest $\b>0$ such that $\cM$ is
completely $\b$-bounded.
\end{Thm}
\begin{proof}
Set $2\b K= T\log\D$ as in Cor. \ref{T}. Note first that
$\cK=\overline{\cM_{sa}\O}$ is equal to $\{\x\in\cH: S\x=\x\}$ where
$S=J\D^{\frac{1}{2}}$ is the Tomita operator. As $T$ commutes both
with $\D$ and $J$, the same is true for $T^{\frac{1}{2}}$ and thus
$T^{\frac{1}{2}}$ commutes with $S$. It follows that
\[
T^{\frac{1}{2}} \cK\subset\cK\ .
\]
We then have for all $\x\in\cK$
\[
2\b(K\x,\x)=(T\log\D\x,\x)=
(\log\D T^{\frac{1}{2}}\x,T^{\frac{1}{2}}\x)\leq 0 .
\]
Clearly the same is true if we replace the von Neumann algebra
$\cM$ by $\otimes_{\mathbb Z} \cM$ (infinite tensor product with respect
to the constant sequence of vectors $\O_n\equiv\O$), the vector $\O$
by $\otimes_{\mathbb Z} \O_n$ and $U(t)$ by $\otimes_{\mathbb Z} U(t)$.

The permutation shift then acts in a strongly cluster fashion on the
latter system.  By the Pusz-Woronowicz theorem \ref{PW} either $K=0$
or this latter system satisfies the KMS condition at some inverse
temperature $\b>0$.  Clearly the KMS condition then holds true also
for the original system.  To prove the last assertion we have to show
that, given $\b>1$, $\cM$ is not completely $\b$-bounded with respect
to $-\log \Delta$ and $\O$.  Indeed, if this were not the case, we
would have $\D^{\b}\leq \text{max}(1,\D)$, by Theorem \ref{PH}, which
is not possible if $\b>1$ unless $\D=1$ in which case also $K=0$.
\end{proof}

Let $\gA$ be a C$^*$-algebra, $\a$ a one-parameter automorphism group
and $\o$ an $\a$-invariant state.  At this point it is natural to say
that $\o$ is completely $\b$-holomorphic if state
$\o\otimes\cdots\otimes\o$ of the $n$-fold (spatial) tensor product
$\gA\otimes\cdots\otimes\gA$ is $\b$-holomorphic with constant $C=1$
for all $n\in\mathbb N$.  We then have:
\begin{Thm}\label{cb=KMS}
Let $\gA$ be a C$^*$-algebra, $\a$ a non-trivial one-parameter
automorphism group and $\o$ an $\a$-invariant state.  The following
are equivalent:
\begin{itemize}
\item[$(i)$] $\o$ is completely $\b$ holomorphic;
\item[$(ii)$] $\o$ satisfies the KMS condition at inverse temperature
$\b_{\rm max}=\sup\{\b>0:\o\, \text{is completely $\b$-holomorphic}\}$
($\b_{\rm max}=+\infty$ means that $\o$ is a ground state).
\end{itemize}
\end{Thm}
\begin{proof}
$(ii)\Rightarrow (i):$ If $\o$ satisfies the KMS condition at inverse
temperature $\b>0$ then it is $\b$-holomorphic with constant $C=1$. It
is also immediate by the inequality (\ref{leq2}) that $\o$ is not
completely $\b'$-holomorphic if $\b'>\b$. For the same
reason it is not completely $\b'$-holomorphic if $\b'>\b$.

If $\o$ is a ground state, then $\o$ is obviously completely $
\b$-holomorphic for all $\b>0$.

$(i)\Rightarrow (ii):$ By considering the GNS representation of $\o$,
it is sufficient to show that Theorem \ref{cb2} holds true without
assuming that $\O$ is separating, but allowing $\o$ to be a ground
state.

Indeed let $E\in\cM$ be the projection onto $\cH_0=\overline{\cM'\O}$.
Clearly $E\cM E$ is $\b$-bounded on $\cH_0$ with respect to $K|_{\cH_0}$
and $\O$. Hence Ad$U(t)|_{\cH_0}$ implements the rescaled modular
group of $(E\cM E, \O)$. Thus $\b' KE=-\log \D$, for some
$\b'\geq\b$, where $\D$ is the modular operator of $E\cM E$ acting on
$\cH_0$.

Moreover, by Corollary \ref{ineq}, we have $e^{-\b K}\leq 1+\D E$.
In particular $e^{-\b K}(1-E)\leq (1- E)$, thus $K(1-E)\geq 0$.

Thus $\b'K=-\log\D E + L$, where $L\equiv \b'K(1-E)$ is positive.  By
the completely $\b$-holomorphic assumption, the same relation holds
true by replacing the system with its tensor product by itself.  This
is possibly only in two cases: either $\log\D=0$, thus $K>0$ and $\o$ is
a ground state, or $L=0$, namely $\O$ is cyclic and separating and
$\b K=-\log\D$.
\end{proof}
\subsubsection{Appendix. The domain of the analytic generator}
We collect here a few properties for selfadjoint operators needed in
text.
\begin{Lemma}\label{aLemma}Let $K$ be a selfadjoint operator on the
Hilbert space $\cH$, $\b>0$, and $\xi$ a vector in $\cH$.
The following are equivalent:
\begin{itemize}
\item[$(i)$] $\xi\in D(e^{-\frac{\b}{2} K})$;
\item[$(ii)$] The function $t\in\mathbb R\to (e^{itK}\xi,\xi)$ extends to a
function
continuous in $\bar S_{\b}$ and holomorphic in $S_\b$;
\item[$(iii)$] For all $\eta\in D(e^{-\frac{\b}{2} K})$, the function
$t\in\mathbb R\to (e^{itK}\xi,\eta)$ extends to a function
in $A(S_\b)$.
\end{itemize}
In this case
$\underset{t\to i\b}{\rm anal.cont.}(e^{itK}\xi,\xi)=\|e^{-\frac{\b}{2}
K}\xi\|^2$.
\end{Lemma}
\begin{proof} Let $E(\l)$ be the family of projections associated
with $K$ by the spectral theorem, namely $K=\int \l \text{d}E(\l)$.
Then $\xi\in D(e^{-\frac{\b}{2} K})$ if and only if
$\int e^{-\l\b} \text{d}\|E(\l)\xi\|^2 < \infty$, i.e. $e^{-\l\b} \in
L^1(\m)$, where $\m(V)=\int_{V}\text{d}(E(\l)\xi,\xi)$ is the
finite Borel spectral measure associated with $\xi$. By the next Sublemma
(with a change of sign of $\beta$) this holds iff $t\to
\hat\m(-t)\equiv\int e^{it\l}
\text{d}\mu(\l)$ is the boundary value of a function holomorphic in $S_\b$ and
continuous in $\bar S_\b$   therefore
$(i)\Leftrightarrow (ii)$.

If $(ii)$ holds, then $t\to \hat\m(-t)$ extends to a function in
$A(S_\b)$, namely $\hat\m(-t)$ is bounded in the strip $S_\b$, because
\[
|\int e^{iz\l}\text{d}\mu(\l)|\leq
\int e^{-\Im z \l}\text{d}\mu(\l)\leq \m([0,\infty))
+\int e^{-\b\l}\text{d}\mu(\l)
\]
for all $z\in \bar S_\b$.

If $\xi,\eta\in D(e^{-\frac{\b}{2} K})$, the function
\[
f(z)\equiv (e^{itK}e^{-sK}\xi,e^{-sK}\eta), \quad z =t+is\in\bar S_\b,
\]
can be checked to belong to $A(S_\b)$ by standard methods, thus
$(ii)\Leftrightarrow (iii)$.
Last assertion follows by the spectral theorem.
\end{proof}
\begin{Sublemma} Let $\mu$ be a finite Borel measure on $\mathbb
R$, $\hat\mu$ its Fourier transform and $\b\in\br$.
The function $t\in\mathbb R\to e^{\b t}$  belongs to
$L^1(\mu)$ if and only if $\hat\mu$ is the boundary value of
a function holomorphic in $S_\b$ and continuous in $\bar S_\b$.

This function is automatically bounded, namely it
belongs to $A(S_{\b})$.
\end{Sublemma}
\begin{proof} If $\l\in\mathbb R\to e^{\b \l}$ belongs to $L^{1}(\m)$,
then also $\l\in\mathbb R\to e^{-iz\l}$ belongs to $L^{1}(\m)$
for all $0\leq \Im z\leq \b$
and  $\hat\mu(z)=\int e^{-iz\l}\text{d}\mu(\l)$ defines
a function in the strip $\bar S_\b$, that can be easily seen to
belong to $A(S_\b)$.

Conversely suppose $\hat\mu$ to be the boundary value
of a function holomorphic in $S_{\b}$ and continuous in $\bar S_\b$.
Decompose $\m$ as $\m_{+}+\m_{-}$,
where the first term is supported in the positive axis, the second
in the negative axis. We have $e^{\b t}\in L^1(\mathbb R,\mu_{-})$ for
any positive $\b$, hence $\hat\m_{-}$ is holomorphic in the upper half
plane.  Therefore we may restrict to the case where $\m$ is supported
in the positive axis.

In this case $\hat\m$ is holomorphic in the lower half plane and in
the strip $S_{\b}$, and is continuous on the real line both from above
and from below, therefore $\hat\m$ extends to a holomorphic function
on $\{\Im z<\b\}$, continuous on the boundary.  Set $\f(x)=\hat\m(ix)$.
Then $\f$ is analytic in $x<\b$, continuous on $x\leq\b$, and for
$x\leq0$ is given by $\f(x)=\int e^{x\l}\text{d}\m(\l)$.  If $x<0$, the
dominated
convergence theorems entails $\f^{(n)}(x)=\int \l^{n}e^{x\l}\text{d}\m(\l)$,
hence  by
monotone convergence we obtain $\f^{(n)}(0)=\int \l^{n}\text{d}\m(\l)$.  The
analyticity implies that for $0\leq x<\b$
\begin{equation*}
\f(x)= \sum_{n=0}^{\infty}\frac{x^{n}}{n!}\int \l^{n}\text{d}\m(\l)
=\int e^{x\l}\text{d}\m(\l)
\end{equation*}
where the last equality follows by monotone convergence.  Again by
monotone convergence and the continuity of $\f$ on the boundary we get
$\f(\b)=\int e^{\b\l}\text{d}\m(\l)$, namely $e^{\b t}\in L^1(\m)$.
\end{proof}
We shall also need the following proposition.
\begin{Prop}\label{U-core} Let $U(t)=e^{itK}$ be a one-parameter
unitary group on a Hilbert space $\cH$ and $\f:\mathbb R\to\mathbb C$
a locally bounded Borel function. If $\cD\subset D(\f(K))$ is a dense,
$U$-invariant linear space, then $\cD$ is a core for $\f(K)$.
\end{Prop}
\begin{proof}
By replacing $\f$ with $|\f|$, we may assume that $\f$ is non-negative.
Let $\xi\in\cH$ be a vector orthogonal to $(\f(K) + 1)\cD$. We have to
show that $\xi=0$. If $f$ is a function in the Schwartz space
$S(\mathbb R)$, we have
\begin{equation}\label{orth}
((\f(K) + 1)f(K)\eta,\xi)=
\int\tilde f(t)((\f(K) + 1)e^{-itK}\eta,\xi)\text{d}t =0 \ ,
\end{equation}
for all $\eta\in \cD$,
where $\tilde f$ the Fourier anti-transform of $f$.

If $f$ is a bounded Borel function with compact support, we may choose
a sequence of smooth functions $f_n$ with compact support such that
$f_n(K)\to f(K)$ weakly, thus eq. (\ref{orth}) holds for such an $f$.

If $g$ is a bounded Borel function with compact support, we may write
$g(\lambda)=(\f(\lambda) + 1)(\f(\lambda) + 1)^{-1}g(\lambda)$,
therefore
\[
(g(K)\eta,\xi)= 0\ .
\]
We can then choose a sequence $g_n$ of such functions such
that $g_n(K)\to 1$ strongly. It follows that $(\eta,\xi)= 0$ for all
$\eta\in \cD$, hence $\xi=0$ because $\cD$ is dense.
\end{proof}
\subsection{Case where further symmetries are present}
We now examine the case where the one-parameter unitary group $\a$
considered before extends to a unitary representation of the
``$ax+b$'' group, where positive translations implement endomorphisms of
the algebra.
\begin{Prop}\label{CR}
Let $\cM$ be a von Neumann algebra on the Hilbert space
$\cH$, $\O$ a cyclic vector for $\cM$ and $U$ a $\O$-fixing
one-parameter unitary group on $\cH$, with generator $K$, implementing
automorphisms of $\cM$. Assume furthermore that there is a
one-parameter unitary group $T$ on $\cH$ such that
\[
T(a)\cM T(-a)\subset \cM,\quad \forall a\geq 0
\]
and satisfying the commutation relations
\begin{equation}\label{commutation}
U(t)T(a)U(-t)=T(e^t a),\quad a,t\in\mathbb R\ .
\end{equation}
If there is a dense $^*$-subalgebra $\gB$ of $\cM$ with
such that $\gB\O\subset D(e^{-\pi K})$ and either

- $U(t)\gB U(-t)=\gB$ and $T(a)\gB T(-a)\subset\gB$, $t\in\mathbb R,
a\in\mathbb R^+$, or

- $\gB$ is $\pi$-bounded with respect to $K$ and $\O$,

\noindent
then the generator $H$ of $T$ is positive.
\end{Prop}
\begin{proof}
Note first that, as a consequence of the commutation relations
(\ref{commutation}), we have $T(a)\O=\O$ for all $a\in\mathbb R$
\cite{GL3}.
By the the criterion given in Proposition \ref{CL} below, it will
suffice to construct a core $\cD$ for $e^{-\pi K}$ such that
$T(a)\cD\subset\cD$ for all $a\geq 0$.

The set $\cD\equiv\gB\O$ is contained in
the domain of $e^{-\pi K}$, and clearly
\[
T(a)\cD=T(a)\gB T(-a)\O\subset\gB\O=\cD,\quad a\geq 0\ .
\]
If moreover $\gB$ is {\rm Ad}$U$-invariant, then
\[
U(t)\cD=U(t)\gB U(-t)\O=\gB\O=\cD,\quad t\in\mathbb R\ .
\]
We may thus apply Lemma \ref{core} below to conclude that $\cD$ is a
core for $e^{-\pi K}$.

On the other hand, if $\gB$ is $\pi$-bounded, then also $\cM$ is
$\pi$-bounded by Lemma \ref{1.1}, thus we are in the previous case
as $\cM$ is {\rm Ad}$U$-invariant.
\end{proof}
\begin{Cor}\label{Cor:KMS}
In the previous Prop. \ref{CR}, suppose that $\gB$ is $\pi$-bounded and
further that $\O$ is separating for $\cM$ and $\mathbb C\O$
are the only $T$-invariant vectors.

Then $\o\equiv (\cdot\,\O,\O)$ is a
KMS state for $\a=\text{\rm Ad}U$ at inverse temperature $\b=2\pi$. Namely
the modular operator associated with $(\cM,\O)$ is  $\Delta =e^{-2\pi K}$.
\end{Cor}
\begin{proof}
Since $U$ implements automorphisms of
$\cM$, it commutes with the modular operator $\D$
associated with $(\cM,\O)$, thus $U(2\pi t)\D^{it}$ is a
one-parameter group of unitaries. We denote by $L$
its self-adjoint generator.

By the modular theory $U$ also commutes with the modular conjugation
$J$ associated with $(\cM,\O)$, thus $Je^{\pi K}J=e^{-\pi K}$.
We then have,
for all $X\in\cM_1$,
\begin{equation}\label{bound}
\|e^{-\pi L}X\O\|=\|e^{\pi K}\D^{1/2}X\O\|=
\|e^{\pi K}JX^*\O\|=\|Je^{-\pi K}X^*\O\|\leq C\ .
\end{equation}
We shall now show that the above bound holds for all $X$ in the unit
ball of the $^{*}$-algebra $\gC\equiv\cup_a T(a)\cM T(-a)$.

By Borchers theorem \cite{Borc1}, $\D^{it}$ has the same
commutation relations (\ref{commutation}) as $U(t)$ with $T(a)$,
thus $e^{itL}$ commutes with $T(a)$. Therefore if $X\in\cM_1$
\begin{multline}
\|e^{-\pi L}T(a)XT(-a)\O\|=\|e^{-\pi L}T(a)X\O\|\\
=\|T(a)e^{-\pi L}X\O\|=\|e^{-\pi L}X\O\|\leq C,
\end{multline}
namely the $\pi$-boundedness property with respect to $L$ holds for $\gC$.

Now, by the following Lemma \ref{irreducible},
$\gC$ is irreducible on $\cH$, thus $L$ is semi-bounded by Cor.
\ref{pure}. As $JLJ = - L$, $L$ is indeed a bounded operator.

We now follow an argument in \cite{BW}. By the Kadison-Sakai derivation
theorem (cf. \cite{Takesaki}), there exists a
selfadjoint element $h\in \cM$, indeed a minimal positive one, such that
\[
e^{ith}Xe^{-ith}= e^{itL}Xe^{-itL}, \quad X\in\cM\ ,
\]
and indeed $\D^{it}h\D^{-it}=h$ by the canonicity of the minimal
positive choice for $h$. Therefore $\D^{it}h\O=h\O$, $t\in\mathbb R$,
and this implies $T(a)h\O=h\O,\, a\in\mathbb R$ \cite{GL3}, thus
$h\O\in\mathbb C\O$ by the uniqueness of the $T$-invariant vector. As $\O$ is
separating, $h\in\mathbb R^+$, thus $h=0$ as $h$ is minimal.

It follows that $\D^{it}=U(-2\pi t)$ for all $t\in \mathbb R$.
\end{proof}
Note that, by an argument of Driessler, see \cite{L}, the von Neumann
algebra $\cM$ in the above corollary is a $III_1$-factor, unless
dim$\cH\leq 1$.
\subsubsection{Appendix. Spectral and irreducibility properties}
We begin to recall a simple lemma.
\begin{Lemma}{\rm \cite{DL}.}\label{E}
Let $\gC$ be a $^*$-algebra on a Hilbert space $\cH$ with cyclic vector
$\O$ and $E$ be the one-dimensional projection onto $\mathbb C\Omega$.
The $^*$-algebra generated by $\gC$ and $E$ is irreducible.
\end{Lemma}
\begin{proof} Let $X\in B(\cH)$ commute with $\gC$ and $E$. Then
$X\in\gC'$ and $X\O=XE\O=EX\O=\lambda\O$ for some $\lambda\in\mathbb C$.
As $\O$ is separating for $\gC'$, then $X=\lambda$ and this entails
the thesis.
\end{proof}
\begin{Lemma}{\rm \cite{L}}
\label{irreducible}
Let $\gC$ be a $^*$-algebra on a Hilbert space $\cH$, $T$ a $n$-parameter
unitary group, $n\geq 1$, such that $T(x)\gC T(-x)= \gC$ for $x\in\mathbb
R^n$. If the
spectrum of $T$ is asymmetric, namely $\sp(U)\cap
-\sp (U)=\{0\}$ and $\O$ is a vector which is cyclic for
$\gC$ and unique $T$-invariant, then $\gC$ is irreducible.
\end{Lemma}
\begin{proof} Let $\cM$ be the weak closure of $\gC$.  Clearly
$T(x)\cM T(-x)=\cM$, hence, by a theorem of Borchers, \cite{Borc2},
see also \cite{Flo}, $T(x)\in\cM$.

By the mean ergodic theorem, the one dimensional
projection $E$ onto $\mathbb\O$ belongs to von Neumann algebra
generated by $\{T(x), x\in\mathbb R^{n}\}$, hence to $\cM$.

Then $\cM=B(\cH)$ by the Lemma \ref{E}.
\end{proof}
\begin{Lemma}\label{core}
Let $\gB$ be a $^*$-algebra, $\O$ a cyclic vector
for $\gB$ and $U(t)=e^{itK}$ a one-parameter, $\O$-fixing unitary group
implementing automorphisms of $\gB$.
If $\gB\O\subset D(e^{-\b K})$ for some $\b>0$,
then $\gB\O$ is a core for $e^{-\b K}$.
\end{Lemma}
\begin{proof}
Set $\cD=\gB\O$ and apply Proposition \ref{U-core}.
\end{proof}
We now recall the criterion for the positivity of the energy discussed
in \cite{BCL}.
\begin{Prop}\label{CL} {\rm \cite{BCL}}. Let $U$ and $T$ be
one-parameter unitary group on a Hilbert space $\cH$ satisfying the
commutation relations (\ref{commutation}). The following are equivalent:
\begin{itemize}
\item[$(i)$] the generator of $T$ is positive;
\item[$(ii)$] there exists a core $\cD$ for $e^{-\pi K}$ such that
$T(a)\cD\subset\cD$ for some (hence for all) $a>0$, where $K$ is the
generator of $U$.
\end{itemize}
\end{Prop}




\section{Minimality of the Hawking temperature}
We now apply our results in the Quantum Field Theory context.
\subsection{A characterization of the Bisognano-Wichmann property}
In the following we shall consider a Poincar\'e covariant net of
von~Neumann algebras on the Minkowski spacetime in the vacuum
representation, indeed an inclusion preserving map
  $$\cS\to\cA(\cS)$$
associating a von~Neumann algebra acting on a given Hilbert space
$\cH$ with each spacelike cone $\cS$ in the Minkowski spacetime
$\mathbb R^{d+1}$, $d\geq 1$, satisfying the following properties:
\begin{itemize}
\item There exists a representation $U$ of the Poincar\'e group $\Ppo$ such
that
$$
U(g)\cA(\cS)U(g)^*=\cA(g\cS)\, ,\quad g\in\Ppo \ .
$$
\item There exists a unique unit vector $\O$ which is
invariant under the action of the Poincar\'e group. $\O$ is cyclic and
separating for the von Neumann algebras associated with wedge
regions.
\end{itemize}
Here a spacelike cone $\cS$ is the cone generated by double cone and
a point in the interior of its spacelike complement.

Note that we do neither assume $\cA$ to be local nor the positivity
of the energy to hold. This last property will indeed follow by our
boundedness condition.

With each wedge region we associate the one-parameter group $\L_W$
of Lorentz boosts preserving $W$. In this way, denoting by $\D_W$
the Tomita operator associated with the von~Neumann algebra $\cA(W)$
and $\O$, the Bisognano-Wichmann relations \cite{BiWi1} take the form
 \begin{equation}\label{BW}
\D_W^{it}=U(\L_W(-2\pi t))\ ,
\end{equation}
where $W$ is a wedge. Clearly, by Poincar\'e covariance, eq. (\ref{BW})
holds for all wedges if it holds for a particular one.
If $W$ is a wedge we denote by $K_W$ the generator
of the one-parameter unitary group $U(\L(t))$.

We now need the following geometric observation, whose proof is
straightforward:
\begin{Lemma}\label{transl} If $\cS$ is an open convex cone  of
$\mathbb R^n$, the set of its translated $\{\cS+x:x\in\mathbb R^{n}\}$
is directed with respect to inclusion.
\end{Lemma}

\begin{Thm}\label{main} Let $W$ be a wedge, $\cS$ a spacelike cone contained
in $W$ and $\gB$ a weakly dense $^*$-subalgebra of $\cA(\cS)$.
The following are equivalent:
\begin{itemize}
\item[$(i)$] The Bisognano-Wichmann relation $\D_W=e^{-2\pi K_W}$ holds.
\item[$(ii)$] $e^{-\pi K_W}\gB_1\O$ is bounded and the
energy-momentum spectrum lays in the forward light cone $\bar V_{+}$.
\item[$(iii)$] $e^{-\pi K_W}\cA(W)_1\O$ is bounded
\item[$(iv)$] $||e^{-\pi K_W}\cA(W)_1\O||\leq 1$.
\end{itemize}
If moreover the boundary of $\cS$ intersects the edge of $W$ in a
half-line, then in $(ii)$ the $\pi$-boundedness is sufficient, namely
it implies the spectrum condition.
\end{Thm}

\begin{proof} In this proof we drop the subscript $W$ on the operators
associated
with $W$.

$(i)\Rightarrow(iv)$: As $\D^{1/2}=e^{-\pi K}$, then
$e^{-\pi K}X\O=\D^{1/2}X\O=JX^*\O$ for all $X\in\cA(W)$,
which immediately implies $(iv)$.

$(iv)\Rightarrow(iii)$ is obvious.

$(iii)\Rightarrow(ii)$: We only need to show the spectrum condition.
By Poincar\`e covariance it is sufficient to show that the positivity
of the generator of a one-parameter group $T$ of light-like
translations associated with $W$; this satisfies
$T(a)\cA(W)T(-a)\subset\cA(W)$, $a\geq 0$, and the commutation
relations (\ref{commutation}) with $U(t)\equiv e^{itK}$.

Thus the positivity property
follows by the criterion in Proposition \ref{CL}.

$(ii)\Rightarrow(i)$:
Since $U(\L(t))$ implements automorphisms of
$\cA(W)$, it commutes with $\D^{is}$, i.e. $U(\L(2\pi t))\D^{it}$ is a
one-parameter group of unitaries. Denote by $L$ its self-adjoint generator
and note that, since by the Tomita-Takesaki theorem $e^{itL}$ commutes with
the modular conjugation $J$ of $(\cA(W),\O)$, we have for all
$X\in\cA(\cS)_1$,
$$
\|e^{-L/2}X\O\|=\|e^{-\pi K}\D^{1/2}X\O\|=\|e^{-\pi K}X^*\O\|\leq C,
$$
namely $\cA(\cS)$ is $\frac{1}{2}$-bounded with respect to $L$ and $\O$.

Moreover, by Borchers commutation relations \cite{Borc1}
and by Tomita-Takesaki theorem,
$e^{itL}$ commutes with all translations. Therefore, denoting by $T$ the
translation unitary group,
$$
\|e^{-L/2}\cA(\cS + x)_1\O\|=\|e^{-L/2}T(x)\cA(\cS)_1\O\|=
\|e^{-L/2}\cA(\cS)_1\O\|\leq C,
$$
namely $\|e^{-L/2}X\O\|\leq C$ for all $X\in\cup_x\cA(\cS+x)$
with $||X||\leq 1$.

By Lemma \ref{transl} $\cup_x\cA(\cS+x)$ is a $^*$-subalgebra
 of $\cB(\cH)$, which is also translation invariant.
By Lemma \ref{irreducible} $\cup_x\cA(\cS+x)$ is thus irreducible,
thus $L$ is semi-bounded by lemma \ref{pure}.

The rest now follows as in the proof of Corollary \ref{Cor:KMS}.

It remains to show the last assertion. Let's then assume that $\cS$
intersects the edge of $W$ in a half line.

By Lemma \ref{1.1} the boundedness of $e^{-\pi K}\gB_1\O$ implies
that also $e^{-\pi K}\cA(\cS)_1\O$ is bounded.

Let $T(s)$ be the one-parameter unitary group of translations
along the edge of $W$. Clearly $T(s)$ commutes with $K$, therefore
if $||e^{-\pi K}\cA(\cS)_1\O||\leq C$ then
\begin{multline}
||e^{-\pi K}\cA(\cS + s)_1\O||=||e^{-\pi K}T(s)\cA(\cS)_1\O||\\
= ||T(s)e^{-\pi K}\cA(\cS)_1\O||=
||e^{-\pi K}\cA(\cS)_1\O||\leq C\ ,
\end{multline}
namely the $\pi$-boundedness condition hold for $\cup_s \cA(\cS + s)$.
As $\cA(\cS + s)$ is a dense $^*$-algebra of $\cA(W)$, by Lemma
\ref{1.1} $||e^{-\pi K}\cA(W)_1\O||\leq C$, thus we obtain all the
properties in the statement by the above proof.
\end{proof}
Let now $\cM$ be a von Neumann algebra on a
Hilbert space $\cH$ and $\O$ a cyclic and separating vector for $\cM$.
For a vector $\xi\in\cH$ we set $\f_{\xi}\equiv
(\cdot\,\O,\xi)|_{\cM}$.
We shall say that a set $Q\subset\cH$ is $L^1$-metrically nuclear
with respect to $\cM$ if
the set of linear functionals $\{\f_{\xi}\in\cM_*: \xi\in Q\}$ is
a metrically nuclear subset of $\cM_*$.

By using a characterization of the split property of Fidaleo \cite{F},
we obtain the following.
\begin{Cor}
Let $W$ be a wedge region and $\cS$ a
spacelike cone  contained in $W$.
Consider the following properties:
\begin{itemize}
\item[$(i)$] The Bisognano-Wichmann property $\D=e^{-2\pi K_W}$ for $W$;
\item[$(ii)$] The split property for $\cA(\cS)\subset\cA(W)$.
\item[$(iii)$] The set $e^{-\lambda K_W}\cA(\cS)_1\O$ is
compact for every $0<\lambda<\pi$, $L^1$-metrically nuclear with respect
to $\cA(W)$ for $\lambda = 1/2$, and
the diameter of $e^{-\lambda K_W}\cA(\cS)_1\O$ is uniformly bounded for
$0<\lambda<\pi$.
\end{itemize}
Then $(i)\,\&\,(ii)\Leftrightarrow(iii)$.
\end{Cor}
\begin{proof} $(i)\&(ii)\Rightarrow(iii)$: By the split property
the set $\D^{\lambda}\cA(\cS)_1\O$ is compact for $0<\lambda<1/2$
\cite{BDL1} and metrically nuclear for $\lambda = 1/4$ \cite{F}. By
the KMS property $||\D^{\lambda}X\O||\leq 1$ for all $X\in\cA(W)_1$
and $0\leq\lambda\leq \frac{1}{2}$ (see \cite{BDL1});
we omit the suffix $W$ on
$\D$ and $K$.

$(iii)\Rightarrow(i)\&(ii)$: Assume first that
the underlying Hilbert $\cH$ space is separable.
Let $X\in\cA(\cS)_1$ and choose a sequence $\{\xi_n\}$ norm dense in
$\cH_1$; the function
\[
F_n(z)\equiv (e^{izK}X\O,\xi_n)
\]
is bounded and holomorphic in the open strip $S_{\pi}$.
Indeed, by the uniform boundedness assumption,
$|F(z)|\leq C$ for some constant $C>0$ independent of $X$ and $\xi_n$.

As $F_n\in H^{\infty}(S_{\pi})$ the limit
$\lim_{\lambda\to\pi}F_n(t+i\lambda)$ exists except for
$t$ in a set $E_n\subset\mathbb R$ of Lebesgue measure zero.
Choose $t_0\notin \cup_n E_n$; as $||e^{i(t_0 +i\lambda)K}X\O||\leq
C$, the weak limit
$\lim_{\lambda\to\pi}e^{i(t_0 +i\lambda)K}X\O$
exists, thus also the weak limit $\lim_{\lambda\to\pi}e^{-\lambda K}X\O$
exists. By the spectral theorem, this implies $X\O\in D(e^{-\pi K})$ and
$||e^{-\pi K}X\O||\leq C$, namely $\cA(\cS)$ is $\pi$-bounded with
respect to $K$ and $\O$ and this entails $(i)$ by the previous theorem.
If $\cH$ is non-separable, it is sufficient to apply the above
argument to the separable Hilbert subspace generated by $f(K)X\O$ as
$f$ varies in the complex continuous functions on $\mathbb R$
vanishing at infinity.

Then $(ii)$ follows because the $L^1$-metrical
nuclearity of $\D^{1/2}\cA(\cS)_1\O$ implies the split property for
$\cA(\cS)\subset\cA(W)$ \cite{F}.
\end{proof}
\begin{rem}
	If the net $\cA$ is local, then the modular conjugation $J_W$ for
	wedges has a geometric action too \cite{GL2} if Theorem \ref{main}
	holds.
\end{rem}
\subsection{One-dimensional nets}
It is convenient to give explicitly a version of the above results in
the context of nets of von Neumann algebras on the real line.

Let $\cI$ denote the set of bounded open non-empty intervals of
$\mathbb R$. We shall consider a net of von Neumann algebras on
$\mathbb R$, namely an inclusion preserving map
\[
I\in\cI\to\cC(I)
\]
from $\cI$ to the von Neumann algebras on a given Hilbert space
$\cH$. If $E\subset\mathbb R$, we denote by $\gC(E)$ the
C$^*$-algebra generated by $\{\cC(I): I\in\cI, I\subset E\}$ and by
$\cC(E)$ the weak closure of $\gC(E)$.

We shall further assume that there exist two one-parameter unitary
groups $T$ and $U$ implementing translations and dilations, namely
for any $I\in\cI$, $a,t\in\mathbb R$,
\begin{align}
T(a)\cC(I)T(a)^*&=\cC(I+a)\label{trasl}
\\
U(t)\cC(I)U(t)^*&=\cC(e^t I)\label{dilat}
\end{align}
satisfying the commutation relations (\ref{commutation}) and leaving
invariant a unique unit vector $\O$ (the vacuum), which is cyclic and
separating for
$\cC(0,\infty)$.

We denote by $K$ the infinitesimal generator of $U$.

\begin{Prop}\label{1dim}
The following are equivalent:
\begin{itemize}
\item[$(i)$] The Bisognano-Wichmann relations $\D=e^{-2\pi K}$ holds,
where $\D$ is the modular operator associated with $(\cC(0,\infty),\O)$.
\item[$(ii)$] There exists $a>0$ and a dense $^*$-subalgebra $\gB$
of $\cC(0,a)$ such that $e^{-\pi K}\gB_1\O$ is bounded.
\item[$(iii)$] There exists $a>0$ and a dense $^*$-subalgebra $\gB$
of $\cC(a,\infty)$ such that $e^{-\pi K}\gB_1\O$
 is bounded.
\item[$(iv)$] There exists $I\in\cI$ and a dense $^*$-subalgebra $\gB$
of $\cC(I)$ such that $e^{-\pi K}\gB_1\O$ is bounded
 and the generator of $T$ is a positive operator.
\item[$(v)$] $e^{-\pi K}\cC(0,\infty)_1\O$ is a bounded set.
\end{itemize}
\end{Prop}
\begin{proof} The proof is similar to the one of Theorem \ref{main}.
We only notice that in this case condition $(ii)$ refers to a bounded
interval $(0,a)$. This is possible because
 $\|e^{-\pi K}\gB_1\O\| \leq C\Rightarrow \|e^{-\pi K}\cC(0,a)_1\O\|\leq C$
 and,  by scaling the interval, this gives
 \[
 \|e^{-\pi K}\cC(0,e^t a)_1\O\|  =\|e^{it K}e^{-\pi K}\cC(0,a)_1\O\|
 =\|e^{-\pi K}\cC(0,a)_1\O\| \leq C\
 \]
for any $t$, thus the boundedness condition holds for $\gC(0,\infty)$.
\end{proof}
\begin{rem}
	If locality is further assumed in Prop. \ref{1dim}, then the net $\cC$
extends to a
	conformal net on $S^1$ \cite{GLW}.
\end{rem}
\subsection{Globally hyperbolic spacetimes}
We shall now discuss our results for a class of stationary black hole
spacetimes, namely globally hyperbolic spacetimes with bifurcate
Killing horizon.  As we shall see, for nets with a boundedness
property the KMS property is equivalent to the existence of a
translation symmetry for the net on the horizon.

Let $\cV$ be a $d+1$ dimensional globally hyperbolic spacetime with a
bifurcate Killing horizon.  An example is given by the
Schwarzschild-Kruskal manifold.  We denote by $\mathfrak h_+$ and
$\mathfrak h_-$ the two codimension 1 submanifolds that constitute the
horizon $\mathfrak h = \mathfrak h_+\cup \mathfrak h_-$.  We denote by
$\cL$ and $\cR$ the ``left and right wedges''.

Let $\k=\k(\cV)$ be the surface gravity, namely, denoting by $\chi$
the Killing vector field, the equation $\nabla g(\chi,\chi)=-2\k\chi$
on $\mathfrak h$, with $g$ the metric tensor, defines a function $\k$
on $\mathfrak h$, that is actually constant on $\mathfrak h$
\cite{KW}.  If $\cV$ is the Schwarzschild-Kruskal manifold, then
$\k(\cV)=\frac{1}{4m}$, where $m$ is the mass of the black hole.  In
this case $\cR$ is the exterior of the Schwarzschild black hole.

In what follows $\cR$ is the actual spacetime, and $\cV$ is to be
regarded as a completion of $\cR$.

Let $\cA(\cO)$ be the von Neumann algebra on a Hilbert space $\cH$ of
the observables localized in the bounded diamond $\cO\subset\cR$.
$\cR\subset\cV$ is a $\L$-invariant region and we assume that the
Killing flow $\L_t$ of $\cV$ gives rise to a one parameter unitary
group $U(t)=e^{iKt}$ implementing automorphisms $\a_t$ of the quasi-local
$C^*$-algebra $\gA(\cR)$ such that  $\a_t(\cA(\cO))=\cA(\L_t(\cO))$.

We now consider a locally normal $\a$-invariant state $\f$ on
$\gA(\cR)$.  The net $\gA$ is assumed to be already in the GNS
representation of $\f$, hence $\f$ is represented by a cyclic vector
$\x$.  Let's denote by $\cR_a$ the wedge $\cR$ ``shifted by''
$a\in\mathbb R$ along, say, $\mathfrak h_+$ (see \cite{GLRV}).  If
$I=(a,b)$ is a bounded interval contained in $\mathbb R_+$, we set
$$
\cC(I)=\gA(\cR_a)''\cap\gA(\cR_b)', \quad 0<a<b\ .
$$
One obtains in this way a net of von Neumann algebras localized on the
horizon parametrized by the intervals of $(0,\infty)$, where the
Killing automorphism group $\a$ acts covariantly by dilations, cf.
equation (\ref{dilat}).  We
denote by $\gC(0,\infty)$ the $C^*$-algebra generated by all
$\cC(a,b), b>a>0$.

We shall say that a one-parameter unitary group $T$ implements
translations on the horizon if equations (\ref{trasl}) and
(\ref{commutation}) are satisfied.

\begin{Cor} Assume that $\cC(0,a)$ is $\frac{\b}{2}$-bounded w.r.t.
$K$ and $\xi$ for some $\b\geq\b_0$ and $a>0$, where
$\b_0=\frac{2\pi}{\k}$ is the inverse of the Hawking temperature.  The
following are equivalent:
\begin{itemize}
	\item[$(i)$] $\f|_{\gC(0,\infty)}$ is a KMS state at Hawking
	inverse temperature $\b_0$.  	
	\item[$(ii)$] There exists a one-parameter unitary group $T$
	implementing translations on the horizon.
\end{itemize}
In this case the generator of $T$ is positive and the net extends to a
conformal net on the line.
\end{Cor}
\begin{proof}
	$(i)\imply(ii)$. The inclusion $\cR_{1}\subset\cR_{0}$ is half-sisd
	modular, therefore the translations (with positive generator) can be
	constructed as in \cite{Wies1}, cf. also Proposition4A.2 in \cite{GLRV}.
	\par\noindent
	$(ii)\imply(i)$. First we extend the net $\cC$ to all intervals in
	$\br$ setting $\cC(a,b)=T(a)^{*}\cC(0,b-a)T(a)$, $a<0$. Clearly
	$T$ and $U$ act as translations and dilations on the net, therefore
	Proposition \ref{1dim} applies, hence the generator of $T$ is positive.
	\par\noindent
	Conformal invariance then follows by \cite{Wies4}.
\end{proof}
\section{Conclusion}
We have seen that states with the $\b$-boundedness conditions are
indeed thermal equilibrium states in certain Quantum Field Theory
contexts. In general the $\b$-boundedness condition selects particular
non-equilibrium steady states, whose meaning is not completely clear.

One may be tempted to use the $\b$-boundedness condition to define the
``local temperature'' of an observable $X$ as the inverse of
$\sup\{\b>0: X\O\in D(e^{-\frac{\b}{2}K})\}$ with the notations in the
text, cf. \cite{DZ}.

This suggests the following physical interpretation.
If a thermodynamical system $\Sigma$ sits in the background of the
black hole, interacts with heat reservoirs at temperature less than
the Hawking temperature, then the black hole is a predominant heat bath
for $\Sigma$ and the state is a thermal equilibrium state at the
black hole background temperature. In particular, the Hawking temperature is
minimal and one cannot cool the system $\Sigma$ down by letting it
interact with an infinite reservoir at lower temperature. However,
the above discussion relies on the assumption that only
$\b$-holomorphic states enter in the game.

The validity of the above picture relies on a clarification of the
role of the $\b$-holomorphic states, namely how large is their class and
how close they are to equilibrium states.
\medskip

\noindent {\bf Acknowledgments.} We would like to thank C. D' Antoni
for conversations.



\end{document}